%----------------------------------------------------------------------------
% HANPAPER_no     = 82 
% HANPAPER_type   = article
% PUBPAPER_author = Y. Bugeaud, G.-N. Han
% PUBPAPER_title  = A combinatorial proof of the non-vanishing of Hankel determinants of the Thue--Morse sequence
% PUBPAPER_journal= ? 
% PUBPAPER_volume = ? 
% PUBPAPER_year   = 2013
% PUBPAPER_pages  = 16 
% MAKEPAPER_output= ps pdf 
% MAKEPAPER_4up   = "-m30 -b-30"
%----------------------------------------------------------------------------

% <<<

%%%%%%%%%%%%%%%%%%%%%%%%%%%%%%%%%%%%%%%%%%%%%%%%%%%%
% typoref.tex. V : January 18, 2000. 
% Author : Anthony PHAN
% Warning : syntaxe +- LaTeX 
% Sources :
% T. Lachand--Robert, ``La Ma\^\i trise de \TeX'',
% R\'ef\'erences crois\'ees;
% latex.ltx's sources;
% and of course the \TeX book.
%%%%%%%%%%%%%%%%%%%%%%%%%%%%%%%%%%%%%%%%%%%%%%%%%%%%%
%
\catcode`@=11
%
% style (look at the behavior of \item dans \bibitem too,
% and at one ,\  in \re@dreferenceslist)
% Feel free to change: 	\bibn@me (title like ``R\'ef\'erences'')
%			\bibliographym@rk (general style)
%
\def\bibn@me{R\'ef\'erences}
\def\bibliographym@rk{\centerline{{\sc\bibn@me}}
	\sectionmark\section{\ignorespaces}{\unskip\bibn@me}
	\bigbreak\bgroup
	\ifx\ninepoint\undefined\relax\else\ninepoint\fi}
%
% Beware of the \bgroup: it will be closed by \endthebibliography
%
% \refsp@ce is the spacing command that appens between multiple
% references.
%
\let\refsp@ce=\ 
\let\bibleftm@rk=[
\let\bibrightm@rk=]
%
% if you want more space between brackets...
%\let\refsp@ce=\thinspace
%\def\bibleftm@rk{[\thinspace}
%\def\bibrightm@rk{\thinspace]}
%
% frenchy stuff
%
\def\numero{n\raise.82ex\hbox{$\fam0\scriptscriptstyle o$}~\ignorespaces}
%
% new variables
%
\newcount\equationc@unt
\newcount\bibc@unt
\newif\ifref@changes\ref@changesfalse
\newif\ifpageref@changes\ref@changesfalse
\newif\ifbib@changes\bib@changesfalse
\newif\ifref@undefined\ref@undefinedfalse
\newif\ifpageref@undefined\ref@undefinedfalse
\newif\ifbib@undefined\bib@undefinedfalse
\newwrite\@auxout
%
% mark an equation
%
\def\eqnum{\global\advance\equationc@unt by 1%
\edef\lastref{\number\equationc@unt}%
\eqno{(\lastref)}}
%
% One can reference anything, just copy the former macro
% and use it so: \machin \label{truc}
% In machin you would have defined \lastref by some number
% or any text.
%
% References macros
%
% The next macros are the core of \ref and \cite commands.
% Its first argument may be ref, pageref or bib.
%
% It is too tricky to be explained.
% (It is a bit recursive.)
% It allows using \cite or \ref or ...
% with arbitrary many arguments,
% for instance:
% \cite{knuth1,knuth2,ma pomme}
%
% First argument is always ref, pageref or bib.
%
\def\re@dreferences#1#2{{%
	\re@dreferenceslist{#1}#2,\undefined\@@}}
\def\re@dreferenceslist#1#2,#3\@@{\def\next{#2}%
	\expandafter\ifx\csname#1@@\meaning\next\endcsname\relax
	??\immediate\write16
	{Warning, #1-reference "\next" on page \the\pageno\space
	is undefined.}%
	\global\csname#1@undefinedtrue\endcsname
	\else\csname#1@@\meaning\next\endcsname\fi
	\ifx#3\undefined\relax
	\else,\refsp@ce\re@dreferenceslist{#1}#3\@@\fi}
%
% notice that the former ``,\refsp@ce'' will separate
% multiple arguments. But beware of spaces
% while defining a reference or calling for it!
%
% tricky thing: \newlabel has two arguments
% {labelname}{{\lastref}{\pageref}}
% The second argument is read as two arguments
% by \newl@bel. This was necessary to get
% a jobname.aux containing the same syntax
% LaTeX would produce and use.
%
\def\newlabel#1#2{{\def\next{#1}\newl@bel#2}}
\def\newl@bel#1#2{%
	\expandafter\xdef\csname ref@@\meaning\next\endcsname{#1}%
	\expandafter\xdef\csname pageref@@\meaning\next\endcsname{#2}}
\def\label#1{{%
	\toks0={#1}\message{ref(\lastref) \the\toks0,}%
	\ignorespaces\immediate\write\@auxout%
	{\noexpand\newlabel{\the\toks0}{{\lastref}{\the\pageno}}}%
	\def\next{#1}%
	\expandafter\ifx\csname ref@@\meaning\next\endcsname\lastref%
	\else\global\ref@changestrue\fi%
	\newlabel{#1}{{\lastref}{\the\pageno}}}}
\def\ref#1{\re@dreferences{ref}{#1}}
\def\pageref#1{\re@dreferences{pageref}{#1}}
%
% bibliography macros
%
\def\bibcite#1#2{{\def\next{#1}%
	\expandafter\xdef\csname bib@@\meaning\next\endcsname{#2}}}
\def\cite#1{\bibleftm@rk\re@dreferences{bib}{#1}\bibrightm@rk}
%
% The argument of \beginthebibliography
% is any sequence of numerals which will represent
% the maximum \item's length. If you have less than 9
% \bibitem's, this argument may be {any numeral}.
% if you have between 100 and 999 \bibitem's
% this argument may be {any three numerals},
% and so on.
%
\def\beginthebibliography#1{\bibliographym@rk
	\setbox0\hbox{\bibleftm@rk#1\bibrightm@rk\enspace}
	\parindent=\wd0
	\global\bibc@unt=0
	\def\bibitem##1{\global\advance\bibc@unt by 1
		\edef\lastref{\number\bibc@unt}
		{\toks0={##1}
		\message{bib[\lastref] \the\toks0,}%
		\immediate\write\@auxout
		{\noexpand\bibcite{\the\toks0}{\lastref}}}
		\def\next{##1}%
		\expandafter\ifx
		\csname bib@@\meaning\next\endcsname\lastref
		\else\global\bib@changestrue\fi%
		\bibcite{##1}{\lastref}
		\medbreak
		\item{\hfill\bibleftm@rk\lastref\bibrightm@rk}%
		}
	}
\def\endthebibliography{\egroup\par}
%
% THE NEXT MACRO MUST BE INCLUDED
% IN THE \BYE COMMAND. FOR INSTANCE:
%
% \catcode`@=11
% \outer\def\bye{\@closeaux
% 	\par\vfill\supereject\end}
% \catcode`@=12
%
\def\@closeaux{\closeout\@auxout
	\ifref@changes\immediate\write16
	{Warning, changes in references.}\fi
	\ifpageref@changes\immediate\write16
	{Warning, changes in page references.}\fi
	\ifbib@changes\immediate\write16
	{Warning, changes in bibliography.}\fi
	\ifref@undefined\immediate\write16
	{Warning, references undefined.}\fi
	\ifpageref@undefined\immediate\write16
	{Warning, page references undefined.}\fi
	\ifbib@undefined\immediate\write16
	{Warning, citations undefined.}\fi}
%
% initialization of jobname.aux
%
\immediate\openin\@auxout=\jobname.aux
\ifeof\@auxout \immediate\write16
  {Creating file \jobname.aux}
\immediate\closein\@auxout
\immediate\openout\@auxout=\jobname.aux
\immediate\write\@auxout {\relax}%
\immediate\closeout\@auxout
\else\immediate\closein\@auxout\fi
%
% Let's read this file and open it out
%
\input\jobname.aux
\immediate\openout\@auxout=\jobname.aux
% this file will be closed by \bye.
%
% That's all, folks!
%

\def\bibn@me{R\'ef\'erences bibliographiques}
%\input typpo
%
%\catcode`@=11
\def\bibliographym@rk{\bgroup}
%
% \bye est modifie pour la biblio et la table des matieres
%
\outer\def\bye{ 	\par\vfill\supereject\end}

\def\Q{{\bf {Q}}}

\overfullrule=0pt

\magnification=1200

  \def\pro{\noindent {\bf{Proof : }}}

\def\house#1{\setbox1=\hbox{$\,#1\,$}%
\dimen1=\ht1 \advance\dimen1 by 2pt \dimen2=\dp1 \advance\dimen2 by 2pt
\setbox1=\hbox{\vrule height\dimen1 depth\dimen2\box1\vrule}%
\setbox1=\vbox{\hrule\box1}%
\advance\dimen1 by .4pt \ht1=\dimen1
\advance\dimen2 by .4pt \dp1=\dimen2 \box1\relax}

\def\sm{\smallskip} \def\ens{\enspace} 

\def\build#1_#2^#3{\mathrel{\mathop{\kern 0pt#1}\limits_{#2}^{#3}}}

\def\date {le\ {\the\day}\ \ifcase\month\or 
janvier\or fevrier\or mars\or avril\or mai\or juin\or juillet\or
ao\^ut\or septembre\or octobre\or novembre\or 
d\'ecembre\fi\ {\oldstyle\the\year}}

\font\fivegoth=eufm5 \font\sevengoth=eufm7 \font\tengoth=eufm10

\newfam\gothfam \scriptscriptfont\gothfam=\fivegoth
\textfont\gothfam=\tengoth \scriptfont\gothfam=\sevengoth

\def\cqfd{\unskip\kern 6pt\penalty 500 \raise 0pt\hbox{\vrule\vbox 
to6pt{\hrule width 6pt \vfill\hrule}\vrule}\par}

\def\pro{\noindent {\it Proof. }}

\def\smallsquare{\vbox{\hrule\hbox{\vrule height 1 ex\kern 1 ex\vrule}\hrule}}
%\def\cqfd{\hfill \smallsquare\vskip 3mm}

%%%%%%%%%%%%%%%%%%%%%%%%%%%%%%%%%%%%%%%%%%%%%%%%%%%%%%%%%%%%%%%%%%%

\vskip 5mm

%>>>

\def\Sym{\hbox{\rm Sym}}
\def\odd{{ odd}}
\def\even{{ even}}
\def\qed{\quad\raise -2pt\hbox{\vrule\vbox to 10pt{\hrule width 4pt
\vfill\hrule}\vrule}}
\def\FourDet#1#2#3#4{{\Bigl|\matrix{#1&#2\cr #3&#4\cr}\Bigr|}}
\def\cqfd{\qed}

\def\pointir{\discretionary{.}{}{.\kern.35em---\kern.7em}\nobreak
\hskip 0em plus .3em minus .4em }

\rightline{May 20, 2014}
\bigskip
\centerline{\bf A combinatorial proof of the non-vanishing of Hankel determinants }

\sm
\centerline{\bf of the Thue--Morse sequence}

\vskip 8mm

\centerline{Yann B{\sevenrm UGEAUD} and Guo-Niu H{\sevenrm AN} \footnote{}{\rm
}}

{\narrower\narrower
\vskip 12mm

\proclaim Abstract. {
In 1998, Allouche, Peyri\`ere, Wen and Wen established that the 
Hankel determinants associated with the Thue--Morse sequence on $\{-1, 1\}$
are always nonzero. Their proof depends on a set of sixteen recurrence relations.
We present an alternative, purely combinatorial proof of the same result.
We also re-prove a recent result of Coons on the non-vanishing of the
Hankel determinants associated to two other classical integer sequences.
}

}

\vskip 6mm

\vskip 5mm

\centerline{\bf 1. Introduction} 

\vskip 6mm

Let ${\cal C} (z)$ be a power series in one variable with
rational coefficients,
$$
{\cal C} (z) = \sum_{k \ge 0}  \, c_k z^k, \quad c_k \in \Q.
$$
For $k \ge 1$ and $p \ge 0$, let
$$
H_k^p ( {\cal C} ) := \left|
\matrix{ c_p & c_{p+1} & \ldots & c_{p+k-1} \cr
c_{p+1} & c_{p+2} & \ldots & c_{p+k} \cr
\ \vdots \hfill & \ \vdots \hfill & \ddots &
\ \vdots \hfill \cr
c_{p+k-1} & c_{p+k} & \ldots & c_{p+2k-2} \cr} \right|
$$
be the $(p, k)$-order Hankel determinant 
associated to ${\cal C} (z)$. For simplicity, we write $H_k ({\cal C})$
instead of $H_k^0 ({\cal C} )$. 
The study of the non-vanishing of $(p,k)$-order Hankel determinants is an interesting
question on its own, and this is the purpose of the present paper for the sequence   %%y
$(c_k)_{k \ge 0}$ being the Thue--Morse sequence. 
However, we start the introduction by pointing out a motivation
coming from Diophantine approximation and concerning the study of rational approximation
to the real numbers ${\cal C} (1/b)$, where $b \ge 2$ is an integer. 

Let $\xi$ be an irrational, real number.
The irrationality exponent 
$\mu(\xi)$ of $\xi$ is the supremum
of the real numbers $\mu$ such that the inequality
$$
\biggl| \xi - {p \over q} \biggr| < {1 \over q^{\mu}}
$$
has infinitely many solutions in rational numbers $p/q$.
It follows from the theory of continued fractions
that $\mu (\xi)$ is always greater than or equal
to $2$, and an easy covering argument shows that 
$\mu (\xi)$ is equal to $2$ for almost
all real numbers $\xi$ (with respect
to the Lebesgue measure). 
Furthermore, Roth's theorem
asserts that the irrationality exponent of every
algebraic irrational number is equal to $2$. The reader is directed
to the monograph \cite{SchmLN} for proofs and refinements of these assertions. 
It is in general a very difficult problem to determine the
irrationality exponent of a given transcendental real number $\xi$, unless
$\xi$ is given by its continued fraction expansion. 
Apart from 
more or less {\it ad hoc} constructions, 
there are only very few examples 
of transcendental numbers $\xi$ whose irrationality
exponent is known. When they
can be applied,
the current techniques allow us most often only to 
get an upper bound for $\mu (\xi)$.

%In particular when $\xi$ is defined by
%its expansion in some integer base $b \ge 2$, we do not 
%generally get enough information to determine the
%exact value of $\mu (\xi)$. Write
%$$
%\xi = \lfloor \xi \rfloor +
%\sum_{k \ge 1} \, {a_k \over b^k},   
%$$ 
%where $a_k \in \{0, 1, \ldots , b-1\}$ for $k \ge 1$
%and $\lfloor \cdot \rfloor$ denotes the integer part function.
Recently, Bugeaud \cite{Bu11} developed a method for 
computing $\mu (\xi)$ when $\xi$ is a  Thue--Morse--Mahler number. 
Let 
$$
{\bf t} = t_0 t_1 t_2 \ldots
%= 0110100110010110100101100110100110010110 \ldots
= 1 \ -1 \ -1 \ 1 -1 \ 1 \ 1 \ -1 \ -1 \ 1 \ 1 \ -1 \ldots 
$$ 
denote the Thue--Morse word on $\{-1, 1\}$ defined by $t_0 = 1$,
$t_{2k} = t_k$ and $t_{2k + 1} =  - t_k$ for $k \ge 0$. 
Alternatively, $t_k = 1$ (resp. $-1$)
if the number of $1$'s in the binary expansion of $k$ is
even (resp. is odd). Let 
$$
{\cal T} (z) := \sum_{k \ge 0} \, t_k z^k
$$
be the generating function of $(t_k)_{k \ge 0}$. It is proved in \cite{Bu11}
that, for every integer $b \ge 2$, the irrationality exponent  of the real number
%(which is transcendental, by an old result of Mahler \cite{Mah29})
$$
%\xi_{{\bf t}, b} = 
{\cal T} (1/b) = 
\sum_{k \ge 0} \, {t_{k} \over b^{k}}
=  1 - {1 \over b} - {1 \over b^2} + {1 \over b^3}  - {1 \over b^4} +
{1 \over b^{5}} + \ldots 
$$ 
is equal to $2$. 
There are two main ingredients in the proof.

A first one is the fact that ${\cal T} (z)$ satisfies
a functional equation, namely
$$
{\cal T} (z) = (1 - z) {\cal T} (z^2) = \prod_{n \ge 0} \bigl( 1 - z^{2^n} \bigr),
$$
a key tool in Mahler's proof \cite{Mah29} that ${\cal T} (1/b)$ is transcendental. 

A second one is the non-vanishing of Hankel determinants
associated with ${\bf t}$, a result established by 
Allouche, Peyri\`ere, Wen and Wen \cite{APWW98}.

\proclaim Theorem APWW.
%Let $\Xi(z)$ be the generating function of the Thue--Morse
%sequence on $\{-1, 1\}$ starting with $1$. 
For every positive integer $k$,
the Hankel determinant $H_k ( {\cal T} )$ is non-zero.

The proof given in \cite{APWW98} is long and difficult. It depends on
a set of sixteen recurrence relations involving the $(p, k)$-order Hankel determinants
and gives additional results on the values 
of the Hankel determinants $H^p_k ( {\cal T} )$.

Subsequently, Coons \cite{Co13} considered the functions 
$$
{\cal F}(z):=\sum_{k \ge 1} f_k z^k=\sum_{n=0}^\infty {z^{2^n} \over 1+z^{2^n}}
\quad\hbox{and}\quad  
{\cal G}(z):=\sum_{k \ge 1} g_k z^k=\sum_{n=0}^\infty {z^{2^n} \over 1-z^{2^n}}.
\eqno (1.1)
$$
He proved that, for every integer $b \ge 2$, the irrationality exponent
of $ {\cal F}(1/b)$ and $ {\cal G}(1/b)$ is equal to $2$. 
To this end, he followed the method of \cite{Bu11}, replacing the use of Theorem APWW
by that of the next result (Theorem 2 of \cite{Co13}).

\proclaim Theorem C.
For every positive integer $k$,
the Hankel determinants $H_k^1 ({\cal F})$ 
and $H_k^1 ({\cal G})$ are nonzero.

Coons' proof of Theorem C is of the same level of difficulty as the one of Theorem APWW. 
It is long and hard to follow. 

The aim of this note is to provide a unified, combinatorial proof
of both Theorems APWW and C. We believe that our approach is much simpler
than that of \cite{APWW98,Co13}. 

Our paper is organized as follows. The key combinatorial result, namely Theorem J, is stated
in Section 2, along with three equivalent lemmas. Complete proofs of these lemmas and theorem
are given in Section 3. We gather in Section 4 some additional statements, which follow
from our approach. Then, in Section 5, we show how Theorems APWW and C 
can be easily derived from Theorem J. 
%The final section is devoted to concluding remarks. 

When nothing else is specified, the notation $a \equiv b$ means that the integers $a$
and $b$ are congruent modulo $2$.

%>>>

\vskip  5mm
\centerline{\bf 2. Permutations and involutions} % <<<
\vskip 5mm

Throughout this text, $N$ denotes the set of non-negative integers. 
We introduce the sets
$$
J= \{(2n+1)2^{2k} -1\mid n,k \in N\} 
= \{ 0, 2, 3, 4, 6, 8, 10, 11,  \ldots \},
$$
and
$$
K=N\setminus J=\{1,5,7,9,13,17,\ldots\} = \{(2n+1)2^{2k+1}-1 \mid n,k \in N\}.
$$

Let $\Sym_m=\Sym_{\{0,1,\ldots, m-1\}}$
be the set of all permutations on $\{0,1,\ldots, m-1\}$.
In this section we prove the following result.

\proclaim Theorem J.
\smallskip
(J1) For every integer $m\geq 1$, the number of 
permutations $\sigma\in\Sym_m$
such that $i+\sigma(i) \in J$ for $i=0,1,\ldots, m-1$,
is an odd number.
\smallskip
(J2) For every integer $m\geq 1$, the number of 
permutations $\sigma\in\Sym_m$
such that $i+\sigma(i) \in J$ for $i=0,1,\ldots, m-2$,
is an odd number.

The proof is based on some combinatorial techniques.
Since we want to enumerate permutations modulo 2, we can delete
suitable pairs of permutations  and the result will not be changed.
The problem is then how to associate a given permutation
with another to form a pair. Two methods are used in the present paper:

(1) taking the inverse $\sigma^{-1}$ of a given permutation $\sigma$; 

(2) exchanging two values by 
letting $\sigma(i):=\sigma(j)$ and $\sigma(j):=\sigma(i)$.

Those two methods are fully described in the proof of Theorem $(J1)$.

\medskip

Throughout this paper we use three representations for permutations: 
the {\it one-line}, {\it two-lines} and the {\it product of disjoint cycles}. 
For example, we write
$$
\sigma\in\Sym_9=516280374=
\pmatrix{012345678\cr 516280374\cr} 
= (0,5)(1)(2,6,3)(4,8)(7). 
$$
We choose to separate the elements of a cycle by commas.

Sometimes we write a list of sets under the two-lines representations. If the
set $A$ is under the index $i$, this means that $i+\sigma(i)\in A$.
For example, the permutations in  $(J1)$ and $(J2)$ are 
$$
\pmatrix{
	0&1&2&3&4&5&6&7&8\cr 
\cdot&\cdot&\cdot&\cdot&\cdot&\cdot&\cdot&\cdot&\cdot\cr
	J&J&J&J&J&J&J&J&J\cr 
} 
$$
and
$$
\pmatrix{
	0&1&2&3&4&5&6&7&8\cr 
\cdot&\cdot&\cdot&\cdot&\cdot&\cdot&\cdot&\cdot&\cdot\cr
	J&J&J&J&J&J&J&J&N\cr 
} 
$$
respectively.
\medskip

An {\it involution} is a permutation $\sigma$ such that $\sigma=\sigma^{-1}$.
In the cycle representation of an involution every cycle is either a
{\it fixed point}
$(b)$ or a {\it transposition} $(c, d)$.
For every set $A$, a transposition $(c, d)$ is said to be an {\it $A$-transposition}
if $c+d\in A$ and $c+d$ is odd. This means also that there is an even number
and an odd number in every $A$-transposition. Generally the order of the
two numbers in a transposition does not matter. However, throughout this paper, 
we always write the {\it even number before the odd number} in 
every $A$-transposition.
%, else some transformation could not be reversible.

\medskip

For finite sets of positive integers $A, B$ and a non-negative integer $f$,
let $\nu(A, f, B)$ be the 
number of involutions in $\Sym_{A}$ such that all transpositions 
are $B$-transpositions and have exactly $f$  fixed points. 
Also let $\nu(A, f/g, B) = \nu(A, f, B) +  \nu(A, g, B)$.
For an infinite set $A$ 
let $A|_m$ be the set composed of the $m$ smallest integers in $A$.

We state below three equivalent lemmas. The first (resp. second, third) assertion of any
of these is equivalent to the first (resp. second, third) assertion of any of the
other two lemmas. 

\proclaim Lemma N.
For $m\geq 1$ we have
$$
\leqalignno{
	\nu(N|_{m}, 0/1, K) &\equiv 1; &(N1)\cr
	\nu(N|_{2m-1}, 1, J) &\equiv 1; &(N2)\cr
	\nu(N|_{2m}, 0/2, J) &\equiv 1. &(N3)\cr
}
$$

Let 
$$
P=\{0,3,4,7,8,11,\ldots\}=\{k \mid k=0,3\pmod 4\}
$$
and 
$$
Q=\{1,2,5,6,9,10,\ldots\}=\{k \mid k=1,2\pmod 4\}.
$$
We define three transformations:
$$
\leqalignno{
\beta &: N\rightarrow N; \qquad k\mapsto 
	\cases{k/2 &if $k$ is even\cr (k-1)/2 &if $k$ is odd\cr} \cr
\gamma &: P\rightarrow N; \qquad k\mapsto 
	\cases{k/4 &if $k$ is even\cr (k-3)/4 &if $k$ is odd\cr} \cr
\delta&: N\rightarrow N; \qquad k\mapsto 
	\cases{k+1 &\quad\ \ if $k$ is even\cr k-1 &\quad\ \ if $k$ is odd\cr} \cr
}
$$

The transformation $\beta$ is extended to the involutions $\sigma$
on $N|_m$ such that all transpositions are $K$-transpositions, by applying $\beta$
on every number in the cycle representation of~$\sigma$. The transformation
$\beta$ for involutions is reversible,
even though $\beta$ on $N$ is not reversible.
For example
$$
%\beta(  (7)(0,5)(3,6)(2)(1,8)(4)) = 
\beta(  (7)(0,5)(6,3)(2)(8,1)(4))
= (3)(0,2)(3,1)(1)(4,0)(2).
$$
We do not know a priori whether the fixed point $3$ is obtained from $6$ or from $7$.
We must look at the transposition $(3, 1)$ first. It is obtained 
from the permutation $(6, 3)$ since
we know that an even number is always before an odd number in the transposition.
Thus, we can recover the $K$-transpositions $(6, 3)(0, 5)(8, 1)$. 
All the other numbers are fixed points, so these are $(7)(2)(4)$.
If $(c, d)$ is a $K$-transposition, then 
$$\beta(c)+\beta(d) = {c+d-1\over 2}={(2n+1)2^{2k+1}-1-1\over 2} 
= (2n+1)2^{2k}-1\in J.$$

\medskip
In the same way,
the transformation $\gamma$ is extended to the involutions $\sigma$ on $P|_m$
such that all transpositions are $J$-transpositions, by applying $\gamma$
on every number in the cycle representation of $\sigma$. Again, the transformation
$\gamma$ for involutions is reversible,
even though $\gamma$ on $P|_m$ is not reversible.
For example
$$
\gamma(  (15)(0,11)(12,7)(4)(16,3)(8))
= (3)(0,2)(3,1)(1)(4,0)(2).
$$
If $(c,d)$ is a $J$-transposition and $c,d\in P$, then 
$$\gamma(c)+\gamma(d) = {c+d-3\over 4}={(2n+1)2^{2k}-1-3\over 4} 
= (2n+1)2^{2k-2}-1\in J.$$
\medskip
In fact we can check that the image sets of $\beta$ and $\gamma$ are 
identical, thus the transformation
$\gamma^{-1}\beta$ is well defined. 
The above comments are still valid if $J$ and $K$ are exchanged.
By the bijection $\gamma^{-1}\beta$, the following lemma is equivalent to Lemma N.

\proclaim Lemma P.
For $m\geq 1$ we have
$$
\leqalignno{
	\nu(P|_{m}, 0/1, J) &\equiv 1; &(P1)\cr
	\nu(P|_{2m-1}, 1, K) &\equiv 1; &(P2)\cr
	\nu(P|_{2m}, 0/2, K) &\equiv 1. &(P3)\cr
}
$$

The third transformation $\delta$ is extended to the set of involutions $\sigma$ on $P|_m$
such that all the transpositions are $J$-transpositions, by applying $\delta$
on every number in the cycle representation of $\sigma$. The transformation
$\delta$ for involutions is reversible since $\delta$ for $N$ is reversible.
For example
$$
\delta(  (15)(0,11)(12,7)(4)(16,3)(8))
=  (14)(1,10)(13,6)(5)(17,2)(9)
$$
If $(c, d)$ is a $J$-transposition then $(\delta(c), \delta(d))$ is still
a $J$-transposition since $c+d = \delta(c) + \delta (d)$. 
\medskip
The above comments are still valid if $J$ and $K$ are exchanged.
Using the bijection $\delta$, we see that the following lemma is equivalent to Lemma P.

\proclaim Lemma Q.
For $m\geq 1$ we have
$$
\leqalignno{
	\nu(Q|_{m}, 0/1, J) &\equiv 1; &(Q1)\cr
	\nu(Q|_{2m-1}, 1, K) &\equiv 1; &(Q2)\cr
	\nu(Q|_{2m}, 0/2, K) &\equiv 1. &(Q3)\cr
}
$$

Since Lemmas $N$, $P$ and $Q$ are equivalent, we prove
(P1), (N2), (N3) in Section 3.
%>>>
\vskip 5mm
\centerline{\bf 3. The proofs}
\vskip 5mm

We begin with several comments.
% which are common for all proofs.
The proofs of all the theorems and lemmas are based on induction
on the lengths of the permutations. 
Small values of $m$ can be easily checked by hand.  %%y
The proof of  $(N2)$ uses $(P1)$, the proof of $(P1)$ uses $(J2)$,
and the proof of $(J2)$ uses again $(N2)$. 
This is not a circular reasoning, because
the length of permutations
is smaller than the length of the original permutations. 
For every permutation $\sigma$, we say that $\sigma$ contains a column 
${\odd \choose \even}$ (in the two-line representation)
if there is some odd number $j$ such that $\sigma(j)$ is even.
We make similar sentence pour 
${\even \choose \even}$,
${\even \choose \odd}$
and
${\odd \choose \odd}$.
For short we say that a permutation $\sigma$ is {``\it in (J1)" }
if $\sigma$ satisfies the conditions described in the statement of Theorem (J1),
and that an involution $\sigma$ is {``\it in (P1)"} if $\sigma$
is an involution on the set $P|_m$ having 0 or 1 fixed point, and all
transpositions are $J$-transpositions. 
%So do for all others Theorems and Lemmas.
%
The following basic facts are easy to verify:
\smallskip
(FJ1) The set $J$ contains all even numbers;
\smallskip
(FJ2) If an odd number $x$ is in $J$, then $x\equiv 3\pmod 4$.
\medskip

{\it Proof of $(J1)$}\pointir
We count  the permutations in (J1)  modulo 2.
% and let $\#$ be their number.
If $\sigma$ contains more than two columns ${\odd \choose \odd}$
in the two-line representation,
select the first two such columns
${i_1 \choose j_1}$
and 
${i_2 \choose j_2}$.
We define another permutation $\tau$ obtained from $\sigma$ by exchanging $j_1$ and $j_2$
in the bottom line. This procedure is reversible.
By (FJ1), it is easy to verify that $\tau$ is also a valid permutation in (J1). 
So that we can delete the pair $\sigma$ and $\tau$, and 
there only remain the permutations containing $0$ or $1$ column 
${\odd \choose \odd}$,
in particular, having $0$ or $1$ odd fixed point.
Similarly with an even number,
there only remain permutations containing $0$ or $1$ column 
${\even \choose \even}$. Consequently, the only remaining permutations have
$0,1$ or $2$ fixed points.
Thanks to the bijection $\sigma\mapsto \sigma^{-1}$,
we need only consider the involutions. 
We can check that
 all transpositions are $J$-transpositions.
If $m$ is odd, then the involution contains one fixed point,
and the number of such involutions is $\nu(N|_{m},1,J)\equiv 1$ by Lemma $(N2)$.
If $m$ is even, then the involution may contain $0$ or $2$ fixed points,
and in the case of $2$ fixed points, we can check that there are exactly 
one odd and one even fixed point.
Hence
the number of such involutions is $\nu(N|_{m},0/2,J)\equiv 1$ by Lemma $(N3)$.
\qed

\medskip

{\it Proof of $(N2)$ and $(N3)$}\pointir
By (FJ2) the two numbers in every $J$-transposition are either both in $P$ or 
both in $Q$. This means that no $J$-transposition takes one number in $P$ and another in $Q$. 
The involutions in $(N2)$ and $(N3)$ are of type
$$
\leqalignno{
	P\qquad\quad &\quad\qquad\qquad Q \cr
	(\bullet) (\bullet \bullet) (\bullet \bullet)  (\bullet \bullet)\quad &
	\Big|\quad
(\bullet) (\bullet) (\bullet \bullet) (\bullet \bullet)  (\bullet \bullet) 
(\bullet \bullet)  \cr
}
$$
The two parts composed by numbers from $P$ and from $Q$ are ``independent".
The cardinalities of the two parts and the number of fixed  points in each side
are characterized by~$m$.
If $m=2k+1$ is odd, then every involution $\sigma$ in $(N2)$ has exactly one 
fixed point, in $P|_k$ or $Q|_{k+1}$, according to the parity of $k$. Hence 
$$
\nu(N|_{2k+1}, 1, J)\equiv
\nu(P|_{k}, 0/1, J)
\times \nu(Q|_{k+1}, 0/1, J)\equiv 1.  
$$
The last equality follows from Lemmas $(P1)$ and $(Q1)$.
If $m=2k$ is even and $k$ is odd, then the cardinalities of $P|_k$ and $Q|_k$ 
are odd. So that there is one fixed point in $P|_k$ and one in $Q|_k$.
Hence
$\nu(N|_{2k}, 0, J)=0$ and
$$
\nu(N|_{2k}, 2, J)\equiv
\nu(P|_{k}, 1, J)
\times \nu(Q|_{k}, 1, J) \equiv 1. 
$$
Again, the last equality follows from Lemmas $(P1)$ and $(Q1)$.
If $m=2k$ is even and $k$ is even, then 
the cardinalities of $P|_k$ and $Q|_k$ 
are even. Three situations may occur: 

(i) no fixed point neither in $P|_k$
nor in $Q|_k$; 

(ii) two fixed points in $P|_k$ and no fixed point in $Q|_k$;

(iii) two fixed points in $Q|_k$ and no fixed point in $P|_k$. 

Hence, we get
$$
\leqalignno{
	\nu(N|_{2k}, 0, J)&\equiv \nu(P|_{k}, 0, J) \times \nu(Q|_{k}, 0, J)
	\equiv 1, \cr
	\nu(N|_{2k}, 2, J)&\equiv \nu(P|_{k}, 0, J) \times \nu(Q|_{k}, 2, J) +
	 \nu(P|_{k}, 2, J) \times \nu(Q|_{k}, 0, J)\equiv 0.  \cr
}
$$
The first equality follows from Lemmas $(P1)$ and $(Q1)$.
The second equality is proved by using the bijection $\delta$ described in Section 2.
\qed

\medskip

{\it Proof of (P1)}\pointir
When $m=2k$ is even, 
the image $\tau=\gamma (\sigma)$ by $\gamma$ (described in Section 2) 
of every permutation $\sigma$ in $(P1)$ 
can be identified with 
a permutation $\rho$ on $\{0,1,\ldots, k-1\}$
such that $i+\rho(i)\in J$ for $i=0, 1, \ldots , k-1$. 
For example, taking $\sigma=(0,3)(7,8)(4,11)$,
we have $\tau=\gamma(\sigma) =\gamma ((0,3)(8,7)(4,11)) = (0,0)(2,1)(1,2)$
which can be identified with the permutation on $\{0,1,2\}$ given by 
$$
\rho=\pmatrix{
	0 & 1 & 2 \cr
	0 & 2 & 1 \cr
}.
$$
The number of such permutations is odd by Theorem $(J1)$.
\smallskip
When $m=2k+1$ is odd, every permutation $\sigma$ in 
$(P1)$ has one fixed point. Apply the transformation
$\gamma$ to $\sigma$, then replace the unique singleton $b$ (which is obtained
from the fixed point of $\sigma$) by $(k, b)$ or $(b, k)$ to form a pair. 
There is a unique way to choose between $(k, b)$ and $(b, k)$ such that
$\tau$ can be identified to a permutation $\rho$ on $\{0,1,\ldots, k\}$.
For example take $\sigma=(0,11)(7,8)(3,12)(4)$,
we have $\tau=\gamma(\sigma) =\gamma ((0,11)(8,7)(12,3)(4))
= (0,2)(2,1)(3,0)(1)$,
the single element is $b=1$, we need replace $(1)$ by $(1,3)$ (not $(3,1)$ of course). 
So that $\tau$ is identified to the following permutation on $\{0,1,2, 3\}$
$$
\rho=\pmatrix{
	0 & 1 & 2 & 3 \cr
	2 & 3 & 1 & 0 \cr
	J & N & J & J \cr
}
$$
or its inverse
$$
\rho^{-1}=\pmatrix{
	0 & 1 & 2 & 3 \cr
	3 & 2 & 0 & 1 \cr
	J & J & J & N \cr
}.
$$
We can verify that one of the two permutations $\rho$ and $\rho^{-1}$ is in $(J2)$ 
(look at the letter $N$ at the third line in the above examples).
The number of such permutions is odd by Theorem (J2).\qed
\medskip

{\it Proof of (J2)}\pointir
For every permutation $\sigma$ in (J2), $i+\sigma(i)\in J$ for
$i \leq m-2$. 
If $\sigma(m-1)$ is odd and $\sigma$ contains a ${\odd \choose \odd}$ column,
let ${j \choose \sigma(j)}$ be the first  ${\odd \choose \odd}$ column. 
We then define another permutation $\tau$
obtained from $\sigma$ by exchanging $\sigma(j)$ and $\sigma(m-1)$.
This procedure is reversible. We can delete the pair $\sigma$ and $\tau$.
Similarly we can delete every permutation $\sigma$ such that
$\sigma(m-1)$ is even and $\sigma$ contains an ${\even \choose \even}$ column.
There only remain two types of permutations. If $\sigma(m-1)$ is odd,
then every number under an odd number in the two-line representation is even. 
If $\sigma(m-1)$ is even,
then every number under an even number is odd, as shown here:
$$
\pmatrix{
		 1 & 2 & 3 & 4 & 5 &  & \cdots &| & m-1 \cr
		 e &   & e &   & e &  & \cdots &| & o \cr
}
\hbox{\quad and \quad}
\pmatrix{
	 0 & 1&  2 &3 & 4 &  & \cdots &| & m-1 \cr
	 o &  &  o &  & o &  & \cdots &| & e \cr
}
$$
The letter ``e" and ``o" represent an even and odd numbers respectively.
There are two cases to be considered.
\smallskip
(I) If $m$ is even, then $m-1$ is odd.
The number of odd numbers in $N|_m$
is equal to the number of even numbers in $N|_m$.
\smallskip
(I.1) If $\sigma(m-1)$ is even, then
the permutations are of the type
$$
\pmatrix{
	0& 1& 2& 3& 4& 5& 6 & | & 7 \cr
	o& e& o& e& o& e& o & | & e \cr
	J& J& J& J& J& J& J & | & N \cr
}
$$
which is equivalent to the product 
$$
\pmatrix{
	0& 2& 4& 6  \cr
	o& o& o& o  \cr
	J& J& J& J  \cr
}
\times
\pmatrix{
	1& 3& 5&  7 \cr
	e& e& e&  e \cr
	J& J& J&  N \cr
}.$$
The first factor is equal to $\nu(N|_m, 0, J)$, and the second factor
is equivalent to
$$
\pmatrix{
	1& 3& 5&   \cr
	e& e& e&  e\cr
	J& J& J&   \cr
},$$
which counts $\nu(N|_{m-1}, 1, J)\equiv 1$ by (N2).
\smallskip
(I.2) If $\sigma(m-1)$ is odd, then
there are as many even numbers as odd numbers.
Since the last column  is ${\odd \choose \odd}$,
there is necessarily a unique column ${\even \choose \even}$,
which will be called ``intruder", and marked by the $*$ sign.
The permutations in $(J2)$ are of type
$$
\sum_{*\ even}\pmatrix{
	0& 1& 2^*& 3& 4& 5& 6 & | & 7 \cr
	o& e& e& e& o& e & o& | & o \cr
	J& J& J& J& J& J& J & | & N \cr
}
$$
The column 2 is the ``intruder",
and the letter $J$ under that column can be replaced by $N$.
Exchanging $\sigma(*)$ and $\sigma(m-1)$ yields
$$
\sum_{*\ \even}\pmatrix{
	0& 1& 2^*& 3& 4& 5& 6 & | & 7 \cr
	o& e& o& e& o& e & o& | & e \cr
	J& J& N& J& J& J& J & | & N \cr
}
$$
which is equivalent to the product 
$$
\left( \, 
\sum_{*\ \even} \pmatrix{
	0& 2^*& 4& 6  \cr
	o& o& o& o  \cr
	J& N& J& J  \cr
}
\, \right)
\times
\pmatrix{
	1& 3& 5&  7 \cr
	e& e& e&  e \cr
	J& J& J&  N \cr
}.$$
The first factor is equal to $\nu(N|_m, 2, J)$ and the second one 
is equal to $\nu(N|_{m-1}, 1, J)\equiv 1$ by (N2).
By (I.1) and (I.2)
the total number of permutations in this case is
$$
\nu(N|_m,0,J) \times\nu(N|_{m-1}, 1, J)
+ \nu(N|_m,2,J)\times \nu(N|_{m-1}, 1, J)
\equiv
\nu(N|_m,0/2,J) 
\equiv 1,
$$
where the last equality follows from (N3).
\smallskip
(II) If $m$ is odd, then $m-1$ is even.
The number of even numbers in $N|_m$
is equal to the number of odd numbers in $N|_m$, plus 1.
\smallskip
(II.1)
If $\sigma(m-1)$ is even, then the permutations are of type
$$
\pmatrix{
	0& 1& 2& 3& 4& 5& 6 & 7 & | & 8 \cr
	o& e& o& e& o& e& o & e & | & e \cr
	J& J& J& J& J& J& J & J & | & N \cr
}
$$
which is equivalent to the product 
$$
\pmatrix{
	0& 2& 4& 6  \cr
	o& o& o& o  \cr
	J& J& J& J  \cr
}
\times
\pmatrix{
	1& 3& 5&  7 & 8 \cr
	e& e& e&  e & e\cr
	J& J& J&  J & N\cr
}.$$
The first factor is equal to $\nu(N|_{m-1}, 0, J)$, and the second factor
is equivalent to
$$
\pmatrix{
	1& 3& 5& 7  &  \cr
	e& e& e& e  &  e\cr
	J& J& J& J  & \cr
}.$$
The number of such permutations is equal to $\nu(N|_{m}, 1, J)$, which is
odd, by (N2).
\smallskip
(II.2) If $\sigma(m-1)$ is odd, the permutations are of type
$$
\sum_{*\ even}\pmatrix{
	0& 1& 2^*&3& 4& 5& 6& 7 &| & 8 \cr
	o& e& e& e&  o& e& o&  e &| & o \cr
	J& J& J& J&  J& J& J& J &| & N \cr
}
$$
Exchanging $\sigma(*)$ and $\sigma(m-1)$ yields
$$
\sum_{*\ \even}\pmatrix{
	0& 1& 2^*& 3& 4& 5& 6& 7 & | & 8 \cr
	o& e& o& e& o& e &  o& e & | & e \cr
	J& J& N& J& J& J& J  & J & | & N \cr
}
$$
which is equivalent to the product 
$$
\left( \, 
\sum_{*\ \even} \pmatrix{
	0& 2^*& 4& 6  \cr
	o& o& o& o  \cr
	J& N& J& J  \cr
}
\, \right)
\times
\pmatrix{
	1& 3& 5&  7 & 8 \cr
	e& e& e&  e  & e\cr
	J& J& J&  J & N \cr
}.$$
The first factor 
counts the number of involutions of $N|_{m-1}$ with $2$ fixed points.
Their number is equal to $\nu(N|_{m-1}, 2, J)$. The second factor
is equal to $\nu(N|_{m}, 1, J)\equiv 1$ by (N2).
By (II.1) and (II.2)
the total number of permutations in this case is
$$
\nu(N|_{m-1},0,J)\times \nu(N|_{m}, 1, J)
+ \nu(N|_{m-1},2,J)\times \nu(N|_{m}, 1, J)
\equiv
\nu(N|_{m-1},0/2,J) 
\equiv 1,
$$
where the last equality follows from (N3).\qed

\vskip 5mm
\centerline{\bf 4. Further results}
\vskip 5mm

For subsets $A,B,C, D$ of $N$ let
$
\FourDet ABCD_m
$
denote the number of permutations $\sigma$ in $\Sym_m$ such that 
$$
i+\sigma(i) \in
\cases{
	A,  &if $i, \sigma(i) \in [0, m-2]$ \cr
	B,  &if $i\in [0, m-2]$, $\sigma(i)=m-1$ \cr
	C,  &if $i=m-1, \sigma(i) \in [0, m-2]$ \cr
	D,  &if $i=\sigma(i) =m-1$.  \cr
}
$$

\proclaim Theorem D.
For every $m \ge 1$, we have
$$
\leqalignno{
	\FourDet JJJJ_m &\equiv 1, &(J1)\cr
	\FourDet JNJN_m &\equiv 1, &(J2)\cr
	\FourDet JNNN_m &\equiv m, &(J3)\cr
	\FourDet KKKK_m &\equiv m+1, &(K1)\cr
	\FourDet KNKK_m &\equiv m+1, &(K2)\cr
	\FourDet KNNK_m &\equiv m+1. &(K3)\cr
}
$$

{\it Remark 1}.  The following basic facts are easy to verify:
\smallskip
(FK1) The set $K$ contains no even integer;
\smallskip
(FK2) The set $K$ contains all integers of the form $4n+1$;
\medskip

{\it Remark 2}. Since $J$ contains all even numbers,
taking $D=J$ is equivalent to taking $D=N$;
Since $K$ does not contain any even number,
taking $D=K$ is equivalent to taking $D=\emptyset$.
\smallskip
{\it Remark 3}. Formula (K1) implies that the Hankel determinants of the 
series
$$
{\cal K}(z)={1\over 1-z} + \sum_{k \ge 0} {z^{2^n-1}\over 1- z^{2^n}}
$$
satisfies $H_m ({\cal K}) \equiv m+1$ for $m \ge 1$. 
\smallskip
{\it Remark 4}. Formulas in Theorem D can be used to
derive other formulas,
for example
$$
\FourDet JKKJ_m \equiv m. \leqno{(Ex)}
$$

{\it Proof of (Ex)}\pointir
We have successively
$$
\leqalignno{
	\FourDet JNN{\emptyset}_m 
&\equiv \FourDet JNNN_m - \FourDet JJJJ_{m-1} 
	\equiv m - 1  \hbox{\quad [By $(J3)$ and $(J1)$]}, &(Ex1) \cr
	\FourDet JNJ{\emptyset}_m 
&\equiv \FourDet JNJN_m - \FourDet JJJJ_{m-1} 
\equiv 1 - 1 \equiv 0, &(Ex2) \cr
	\FourDet JJJ{\emptyset}_m 
&\equiv \FourDet JJJN_m - \FourDet JJJJ_{m-1} 
\equiv 1 - 1 \equiv 0, &(Ex3) \cr
	\FourDet JNK{\emptyset}_m 
	&\equiv \FourDet JNN{\emptyset}_m - \FourDet JNJ{\emptyset}_m 
	\equiv (m-1) - 0  \equiv m-1, &(Ex4) \cr
\FourDet JJK{\emptyset}_m 
	&\equiv \FourDet JJN{\emptyset}_m - \FourDet JJJ{\emptyset}_m 
	\equiv 0 - 0  \equiv 0, &(Ex5) \cr
\FourDet JKK{\emptyset}_m 
&\equiv \FourDet JNK{\emptyset}_m -  \FourDet JJK{\emptyset}_m
\equiv (m-1) - 0 \equiv m-1, &(Ex6) \cr
\FourDet JKKJ_m 
&\equiv \FourDet JKK{\emptyset}_m + \FourDet JJJJ_{m-1} 
\equiv (m-1) + 1\equiv m.  \qed&(Ex) \cr
}
$$

The first two identities $(J1)$ and $(J2)$ are simply reproduced from  Theorem J.
\smallskip

{\it Proof of $(J3)$}\pointir 
By the same arguments as in the proof of (J1), the permutations to
enumerate are the involutions $\sigma$ such that

(i) $\sigma$ has $0,1,2$ fixed points;

(ii) every transposition in $\sigma$ is a $J$-transposition, except
the one which contains $m-1$, if it exists;

(iii) If $m-1$ is in a transposition $(m-1, b)$, then $m-1+b$ is odd.
\smallskip

(J3.1) If $m$ is odd, then there is exactly one fixed point.
If $m-1$ is this fixed point, then the number of involutions is
$\nu(N|_{m-1}, 0, J)$. 
If $m-1$ is in a transposition $(b_1, m-1)$
and the fixed point is $(b_2)$, then 
deleting  $m-1$ is a transformation reversible.
We get an involution of $\Sym_{m-1}$ with two fixed points $b_1, b_2$.
The number of such involutions is $\nu(N|_{m-1}, 2, J)$.
Hence the total number of involutions
is $\nu(N|_{m-1}, 0, J) + \nu(N|_{m-1}, 2, J) = \nu(N|_{m-1}, 0/2, J)$, which is an odd number
by Theorem (N3).
\smallskip
(J3.2) If $m$ is even, three situations have to be considered.

(i) If the involution $\sigma$ has no fixed point,
then the number of involutions is $\nu(N|_{m-1}, 1, 0)$.

(ii) If the involution $\sigma$ has $2$ fixed points and $m-1$ is a fixed point, then 
the number of involutions is also $\nu(N|_{m-1}, 1, 0)$.

(iii)
If $m-1$ is in a transposition $(m-1, b_1)$ and  the two fixed points
are $(b_2)$ and $(b_3)$, then $b_2+b_3$ is odd. Exchanging $b_1$ with
$b_2$ (resp. $b_3$) if $b_2\equiv b_1 \pmod2$ (resp. $b_3\equiv b_1\pmod2$)
is a transformation reversible. Hence the total number of
involutions is equal to $\nu(N|_{m-1}, 1, 0)+\nu(N|_{m-1}, 1, 0)$
plus an even number, which is $0$ modulo 2.\qed
\smallskip

{\it Proof of $(K1)$}\pointir 
By (FK1) every permutation in $(K1)$ contains no fixed point.
Since we need only count the involutions, the number of such permutations
is 0 if $m$ is odd. If $m=2k$ is even,
the permutations are of type
$$
\pmatrix{
	0& 1& 2& 3& 4& 5& 6 &   7 \cr
	o& e& o& e& o& e& o &   e \cr
	K& K& K& K& K& K& K &   K \cr
}
$$
which is equivalent to the product  
$$
\pmatrix{
	0& 2& 4& 6  \cr
	o& o& o& o  \cr
	K& K& K& K  \cr
}
\times
\pmatrix{
	1& 3& 5&  7 \cr
	e& e& e&  e \cr
	K& K& K&  K \cr
}.$$
The two factors are equal to $\nu(N|_m, 0, K)\equiv 1$ by (N1).
\qed
\medskip
{\it Proof of $(K2)$}\pointir 
If $m$ is even,
all permutations $\sigma$ in $\FourDet KNKN_m$ are of type
$$
\pmatrix{
	0& 1& 2& 3& 4& 5& 6 &   7 \cr
	o& e& o& e& o& e& o &   e \cr
	K& K& K& K& K& K& K &   N \cr
}
$$
which is equivalent to the product  
$$
\pmatrix{
	0& 2& 4& 6  \cr
	o& o& o& o  \cr
	K& K& K& K  \cr
}
\times
\pmatrix{
	1& 3& 5&  7 \cr
	e& e& e&  e \cr
	K& K& K&  N \cr
}.$$
The number of permutations of the type given by
the first factor is equal to $\nu(N|_{m}, 0, K)\equiv 1$
and that of the type given by the second one
is equal to $\nu(N|_{m-1}, 1, K)\equiv 1$ by (N1).
\smallskip
If $m$ is odd,
all permutations $\sigma$ in $\FourDet KNKN_m$ are of type
$$
\pmatrix{
	0& 1& 2& 3& 4& 5& 6 &   7 &8 \cr
	o& e& o& e& o& e& o &   e &e\cr
	K& K& K& K& K& K& K &   K &N\cr
}
$$
which is equivalent to the product 
$$
\pmatrix{
	0& 2& 4& 6  \cr
	o& o& o& o  \cr
	K& K& K& K  \cr
}
\times
\pmatrix{
	1& 3& 5&  7 & 8\cr
	e& e& e&  e & e\cr
	K& K& K&  K & N\cr
}.$$
The number of permutations of the type given by
the first factor is equal to $\nu(N|_{m-1}, 0, K)\equiv 1$ 
and that of the type given by the second one
is equal to $\nu(N|_{m}, 1, K)\equiv 1$ by (N1).
Hence 
$$
\FourDet KNKN_m\equiv 1
$$
and
$$
\FourDet KNKK_m
\equiv \FourDet KNKN_m - \FourDet KKKK_{m-1}
\equiv 1 - (m-1 +1) \equiv m+1.
\qed
$$

\medskip
{\it Proof of $(K3)$}\pointir 
We need only count the involutions without fixed points, so that the 
number of such permutations is 0 when $m$ is odd.
If $m=2k$ is even, all transpositions are necessarily $K$-transpositions, 
except the one which contains $m-1$. Removing  $m-1$ yields an involution
with one fixed point. The number of such involutions is 
$\nu(N|_{m-1}, 1, K)\equiv 1$ by (N1).\qed

\medskip
Furthermore, in Lemma (N3) we have a mixed formula for $0$ or $2$ fixed points.
From the proof of $(N2)$ and $(N3)$ we can separate it into two more precise 
formulas.

\proclaim Lemma N'.
We have
$$
\leqalignno{
	\nu(N|_{2k}, 0, J) &\equiv k+1, &(N4)\cr
	\nu(N|_{2k}, 2, J) &\equiv k.  &(N5)\cr
}
$$

Using the transfomations $\beta, \gamma, \delta$, we see
that Lemma N' is equivalent to

\proclaim Lemma P'.
We have
$$
\leqalignno{
	\nu(P|_{2k}, 0, K) &\equiv k+1, &(P4)\cr
	\nu(P|_{2k}, 2, K) &\equiv k.  &(P5)\cr
}
$$

\proclaim Lemma Q'.
We have
$$
\leqalignno{
	\nu(Q|_{2k}, 0, K) &\equiv k+1, &(Q4)\cr
	\nu(Q|_{2k}, 2, K) &\equiv k. &(Q5)\cr
}
$$

\vskip 5mm
\centerline{\bf 5. Completion of the combinatorial proofs of Theorems APWW and C}
\vskip 5mm

Recall that 
$$
J= \{(2n+1)2^{2k} -1\mid n,k \in N\} 
= \{ 0, 2, 3, 4, 6, 8, 10, 11,  \ldots \} 
$$
and
$$
K=N\setminus J=\{1,5,7,9,13,17,\ldots\} = \{(2n+1)2^{2k+1}-1 \mid n,k \in N\}.
$$

First, we notice that, for $k \ge 1$, the integer 
$g_k$ is equal to the $2$--adic valuation of $2k$, known sometimes 
as the {\it ruler function}. For $k\ge  1$, we have the recursion
$$
g_{2k+1} =1\qquad\hbox{and}\qquad g_{2k} =1+g_k.
$$ 
This sequence starts as
$$
(g_k)_{k \ge 1}= 1,2,1,3,1,2,1,4,1,2,1,3,1,2,1,5,1,2,1,3,1,2,1,4,\ldots 
$$ 

The proofs of Theorems APWW and C combine Theorem J with the
next two lemmas.

\proclaim Lemma 1.
For $k \ge 0$, the integer $g_{k+1}$ is odd if, and only if, $k$ is in $J$.

\pro
Let $k$ be a non-negative integer.
It follows from (1.1) that $g_{k+1}$ is equal to the number of pairs $(j, m)$ of integers
$j \ge 1$, $m \ge 0$ such that $j 2^m = k + 1$. 
In particular, writing $k + 1 = j_0 2^{m_0}$ with $j_0$ odd, we see that
$$
k + 1 = j_0 2^{m_0} = (2  j_0) 2^{m_0-1} = \ldots = (2^{m_0} j_0 ) 2^0,
$$
showing that $g_{k+1} = m_0 + 1$. Consequently, $g_{k+1}$ is odd if, and only if,
$m_0$ is even, that is, if and only if, $k$ is in $J$. 
\cqfd

\proclaim Lemma 2.
For $k \ge 0$, the integer $\delta_k := (t_{k+1} - t_k)/2$
is odd if, and only if, $k$ is in $J$.

\pro
For $\ell \ge 0$, we have 
$$
\delta_{2 \ell} = (t_{2 \ell +1} - t_{2 \ell})/2 = -t_{\ell} = \pm 1,
$$
thus $\delta_{2 \ell}$ is odd. 
For an odd integer $j \ge 1$ and an integer $k \ge 1$, we have
$$
\delta_{j 2^{2k} - 1} \equiv  (t_{j 2^{2k}} - t_{j 2^{2k} - 1})/2
\equiv  (t_{j 2^{2k-1}} + t_{j 2^{2k-1} - 1})/2  \equiv 
1 + \delta_{j 2^{2k-1} - 1} \ens \pmod{2}
$$
and, likewise,
$$
\delta_{j 2^{2k-1} - 1} \equiv  1 +  \delta_{j 2^{2(k-1)} - 1} \ens \pmod{2}.
$$
Consequently, an immediate induction shows that,
for any odd integer $j$ and any integer $\ell \ge 0$, the integer
$\delta_{j 2^{\ell} - 1}$ is odd if, and only if, $\ell$ is even, that is,
if, and only if, $j 2^{\ell} - 1$ is in $J$. 
This proves the lemma.
\cqfd

\bigskip
\goodbreak

\noindent {\it Proof of Theorem C.}

Let ${\cal J} (z)$ be the power series
$$
{\cal J} (z) = \sum_{k \ge 0}  \, j_k z^k = 1 + z^2 + z^3 + z^4 + z^6 + \ldots
$$
defined by $j_k=1$ if  $k\in J$ and $j_k = 0$ otherwise. 
Let $k$ be a positive integer.
By the definition of determinant, $H_k(\cal J)$ is equal to
$$
\sum_{\sigma\in\Sym_k} (-1)^{inv(\sigma)} j_{0+\sigma(0)}
j_{1+\sigma(1)}\cdots
j_{k-1+\sigma(k-1)},
$$
where $inv(\sigma)$ is the number of inversions of $\sigma$.
The product $j_{0+\sigma(0)} j_{1+\sigma(1)}\cdots j_{k-1+\sigma(k-1)}$
is equal to 1 if $i+\sigma(i)\in J$ for $i = 0, 1, \ldots , k-1$, and is equal to 0 otherwise.
Hence, the above Hankel determinant $H_k(\cal J)$ is equal modulo 2 to
the number of
permutations $\sigma\in\Sym_k$
such that $i+\sigma(i) \in J$ for all $i=0,1,\ldots, k-1$,
which is equal to 1 modulo 2 by Theorem (J1).
We deduce from Lemma 1 that $H_k({\cal J}) \equiv H_k^1 ({\cal G})$.
Consequently, the Hankel determinant 
$H_k^1 ({\cal G})$ is always an odd integer. Furthermore, it follows
from (1.1) that the coefficients of the power series ${\cal F} (z)$ and
${\cal G} (z)$ are congruent modulo 2. This implies
that $H_k^1 ({\cal F}) \equiv  H_k^1 ({\cal G})$, proving that the Hankel determinant 
$ H_k^1 ({\cal F})$ is always an odd integer.   \cqfd

\bigskip

\noindent {\it Proof of Theorem APWW.}

Let $k$ be a positive integer and consider the 
Hankel determinant $H_k ({\cal T})$ of the matrix
$$
\pmatrix{ t_0 & t_{1} & \ldots & t_{k-1} \cr
t_{1} & t_{2} & \ldots & t_{k} \cr
\ \vdots \hfill & \ \vdots \hfill & \ddots &
\ \vdots \hfill \cr
t_{k-1} & t_{k} & \ldots & t_{2k-2} \cr}. 
$$
For $j=1, \ldots , k-1$, subtracting the $(j+1)$-th column from the $j$-th column, we see that
$H_k ({\cal T})$ is equal to the determinant of the matrix
$$
\pmatrix{ t_0 - t_1 & t_{1} - t_2 & \ldots & t_{k-2} - t_{k-1} & t_{k-1} \cr
t_{1} - t_2 & t_{2} - t_3 & \ldots & t_{k-1} - t_k & t_{k} \cr
\ \vdots \hfill & \ \vdots \hfill & \ddots &  \vdots &
\ \vdots \hfill \cr
t_{k-1} - t_k & t_{k} - t_{k+1}& \ldots & t_{2k-3} - t_{2k - 2} & t_{2k-2} \cr}. 
$$
By Lemma 2, this determinant is equal to $2^{k-1}$ times an integer
which has the same parity as the determinant of the matrix
$$
\pmatrix{ j_0 & j_1 & \ldots & j_{k-2} & t_{k-1} \cr
j_1 & j_2 & \ldots & j_{k-1} & t_{k} \cr
\ \vdots \hfill & \ \vdots \hfill & \ddots & \vdots &
\ \vdots \hfill \cr
j_{k-1} & j_{k}& \ldots & j_{2k - 3} & t_{2k-2} \cr},
$$
where $j_0, j_1, \ldots$ are defined in the proof of Theorem C. 
Now, we apply Theorem (J2) to deduce that the determinant of the latter
matrix is an odd integer. It then follows that the 
Hankel determinant $H_k ({\cal T})$ is a nonzero integer. \cqfd

\vskip 8mm

\centerline{\bf References}

\vskip 3mm

\beginthebibliography{999}

\bibitem{APWW98}
J.-P. Allouche, J. Peyri\`ere, Z.-X. Wen, and Z.-Y. Wen, 
{\it Hankel determinants of the Thue--Morse sequence},
Ann. Inst. Fourier (Grenoble)  48  (1998),  1--27.

\bibitem{Bu11} 
Y. Bugeaud.
{\it On the irrationality exponent of the Thue--Morse--Mahler numbers},
Ann. Institut Fourier (Grenoble) 61 (2011), 2065--2076.

%\bibitem{BuVa13}
%P. Bundschuh and K. V\"a\"an\"anen,
%{\it Algebraic independence of the generating functions of SternÕs sequence and of its twist},
%J. Th\'eor. Nombres Bordeaux 25 (2013), 43--57. 

\bibitem{Co13}
M. Coons,
{\it On the rational approximation of the sum of the reciprocals of the Fermat numbers},
The Ramanujan Journal 30 (2013), 39--65.

%\bibitem{CoVr12}
%M. Coons and P. Vrbik,
%{\it An irrationality measure for paperfolding numbers},
%Journal of Integer Sequences 15 (2012), Article 12.1.6, 1--10.

\bibitem{Mah29}
K. Mahler, 
{\it Arithmetische Eigenschaften der L\"osungen 
einer Klasse von Funktionalgleichungen},
Math. Ann. 101 (1929), 342--366. Corrigendum 103 (1930), 532.

\bibitem{SchmLN}
W. M. Schmidt,
Diophantine Approximation.
 Lecture Notes in Math.  {785}, Springer, Berlin, 1980.

\endthebibliography

\bigskip\bigskip
\hbox{\vtop{\halign{#\hfil\cr
Yann Bugeaud\cr
IRMA, UMR 7501\cr
Universit\'e de Strasbourg et CNRS\cr
7, rue Ren\'e Descartes\cr
67084 Strasbourg, France\cr\noalign{\smallskip}
{\tt bugeaud@math.unistra.fr}\cr}}
\qquad
\vtop{\halign{#\hfil\cr
Guo-Niu Han\cr
IRMA, UMR 7501\cr
Universit\'e de Strasbourg et CNRS\cr
7, rue Ren\'e Descartes\cr
67084 Strasbourg, France\cr
\noalign{\smallskip}
{\tt guoniu.han@unistra.fr}\cr}}}

\goodbreak
\bye